# New Lagrangian Relaxation Approach for the Discrete Cost Multicommodity Network Design Problem

Nesrine Bakkar Ennaifer[1,2], Safa Bhar Layeb[1,3], Farah Mansour Zeghal[1,4]

[1]*Université de Tunis El Manar, Ecole Nationale d'Ingénieurs de Tunis*
*UROASIS : Optimisation et Analyse des Systèmes Industriels et de Services*
*BP 37 Le Belvédère, 1002, Tunis, Tunisie*
[2]nesrine.bakkar@gmail.com , [3]bhar_safa@yahoo.fr ,[4]farah.zeghal@enit.rnu.tn

Abstract: We aim to derive effective lower bounds for the Discrete Cost Multicommodity Network Design Problem (DCMNDP). Given an undirected graph, the problem requires installing at most one facility on each edge such that a set of point-to-point commodity flows can be routed and costs are minimized. In the literature, the Lagrangian relaxation is usually applied to an arc-based formulation to derive lower bounds. In this work, we investigate a path-based formulation and we solve its Lagrangian relaxation using several non-differentiable optimization techniques. More precisely, we devised six variants of the deflected subgradient procedures, using various direction-search and step-length strategies. The computational performance of these Lagrangian-based approaches are evaluated and compared on a set of randomly generated instances, and real-world problems.

## 1 INTRODUCTION

The Network Design Problems (NDP) represent an important class of combinatorial optimization problems arising in a wide variety of real-life situations such as telecommunication networks, supply chains and logistic networks, power delivery network planning and aircraft assignment problems, to quote for a few. In this paper, we address a variant of the NDP called the Discrete Cost Multicommodity Network Design Problem (DCMNDP) which is defined as follows. Given a connected undirected graph $G=(V, E)$ where $V$ is a set of $n$ nodes and $E$ is a set of $m$ edges, there is a set of facilities that can be installed on each edge $e$, $e \in E$. Each facility $l$, $l=1...L_e$, $e \in E$, is characterized by a capacity $u_e^l$ and a fixed cost $f_e^l$ that are step-increasing functions, i.e. $u_e^1 < u_e^2 < \ldots < u_e^{L_e}$, and $f_e^1 < f_e^2 < \ldots < f_e^{L_e}$. Each capacity represents the maximum total bidirectional flow that may circulate on the corresponding edge. Moreover, a set of $K$ distinct point-to-point multicommodity flow requirements has to be routed across the network. Each commodity $k$, $k=1...K$, is characterized by its flow value $D^k$, its source node $s^k$ and its sink node $t^k$. The DCMNDP requires designing a minimum-cost network that allows routing all the commodity demand flows between their respective endpoints, while installing at most one facility on each edge. It is noteworthy that the DCMNDP is NP-hard in the strong sense (Johnson et al., 1978).

In the literature, the DCMNDP was first addressed by Minoux (Minoux, 1989). In 1999, Gabrel et.al (Gabrel et al., 1999) developed an exact Bender partitioning procedure and in 2003, they presented several fast greedy heuristics (Gabrel et al., 2003). Later on, Mrad and Haouari have described an exact elaborated constraint generation approach to solve to optimality instances with up to 50 nodes and 100 edges (Mrad & Haouari, 2008). Recently, the variant of the DCMNDP with demand uncertainty was investigated and a Benders decomposition procedure was proposed (Lee et al., 2013).

The special case of the NDP, where single-facility and per unit costs are considered, namely, the multicommodity capacitated fixed-charge network design problem (MCNDP), has drawn the attention of many researchers during the last decades. Indeed, we quote but not exhaustively that valid inequalities (Chouman et al., 2011), exact methods (Hewitt et al., 2012; Bärmann et al., 2013; Gendron & Larose, 2014) and heuristic approaches (Ghamlouche et al., 2003, 2004; Hewitt et al., 2010) were extensively investigated. To derive lower bounds for multicommodity network design problems, the Lagrangian relaxation was successfully used. We refer to the work of Gendron et al. for a survey on relaxations techniques that were applied to the MCNDP (Gendron et al., 1999). Based on this work, Crainic et al. have proposed a calibrating and comparative study of bundle against subgradient methods (Crainic et al., 2001). Recently, extensive computational studies on applying bundle methods (Frangioni & Gorgone, 2014) and subgradient algorithms (Frangioni et al., 2015) were carried. These three last studies were based on an arc-flow formulation of the MCNDP.

The objective of this paper is twofold: First, we aim to derive effective lower bounds for the DCMNDP by applying the Lagrangian relaxation to an arc-path formulation instead of the commonly used arc-node formulation, and second we seek to analyse the computational performance of several deflected subgradient algorithms for solving Lagrangian dual problem.

The reminder of this paper is organized as follows. In section 2, we describe an arc-path formulation of the DCMNDP. In section 3, we present the Lagrangian relaxation of the proposed formulation. Section 4 is dedicated to explore six variants of subgradient algorithms using various direction-search and step-length strategies. In section 5, we report the results of extensive computational experiments carried out on randomly generated instances and real-world problems.

## 2 ARC-PATH FORMULATION

In this section, we describe an arc-path formulation of the DCMNDP. For each commodity $k$, $k=1,...,K$, let $P^k$ be the set of paths originating at the source node $s^k$ and ending at the sink node $t^k$ in the graph $G=(V, E)$. We define the following decision variables:

- $z_r^k$: a nonnegative continuous variable that corresponds to the flow quantity of commodity $k$ that circulates on path $r \in P^k$, $k=1,...,K$,

- $y_e^l$: a binary variable that takes the value $1$ if the facility $l$, $l=1,...,L_e$, is installed on edge $e$, $e \in E$, and $0$ otherwise.

We denote by $a_{er}^k$ the binary constant that takes value 1 if edge $e \in E$ is included in the path $r \in P^k$, of commodity $k$, $k=1, ..., K$, and $0$ otherwise.

Using this notation, an Arc-Path model (**AP**) for the DCMNDP can be stated as follows:

$$(AP): \text{Minimize} \sum_{e \in E} \sum_{l=1}^{L_e} f_e^l y_e^l \qquad (1)$$

subject to:

$$\sum_{l=1}^{L_e} y_e^l \leq 1 \qquad \forall e \in E, \qquad (2)$$

$$\sum_{r=1}^{|P^k|} z_r^k \geq D^k, \quad k=1,...,K, \qquad (3)$$

$$\sum_{k=1}^{K} \sum_{r=1}^{|P^k|} a_{er}^k z_r^k \leq \sum_{l=1}^{L_e} u_e^l y_e^l, \quad \forall e \in E, \qquad (4)$$

$$z_r^k \geq 0, \quad \forall r \in P^k, \quad k=1,...,K, \qquad (5)$$

$$y_e^l \in \{0,1\}, \qquad \forall e \in E, \quad l=1,...,L_e. \qquad (6)$$

The objective (1) is to minimize the total installation cost. Constraints (2) require that at most one facility should be installed on each edge $e$, $e \in E$. Constraints (3) express the total flow requirements for each commodity $k$, $k=1,...,K$. Constraints (4) express the total capacity constraint on each edge $e$, $e \in E$. Constraints (5) and (6) are non-negativity and binary restrictions imposed on $z$ and $y$ variables, respectively.

Obviously, (**AP**) is a linear mixed-integer programming formulation with an exponential number of variables.

## 3 LAGRANGIAN RELAXATION

To perform a Lagrangian relaxation to the arc-path formulation (**AP**), we propose to dualize the capacity constraints (4). For each edge $e$, $e \in E$, let $w_e$ denote the corresponding nonnegative Lagrangian dual variable.

Thus, we obtain the following Lagrangian dual:

$$\theta(w) = \text{Minimize} \sum_{e \in E} \sum_{l=1}^{L_e} \tilde{c}_e^l y_e^l + \sum_{k=1}^{K} \sum_{r=1}^{|P^k|} \tilde{h}_r^k z_r^k \quad (7)$$

subject to: (2), (3), (5), and (6), where

$$\tilde{c}_e^l = f_e^l - w_e u_e^l, \forall e \in E, \quad l = 1,\dots,L_e, \quad (8)$$

$$\tilde{h}_r^k = \sum_{e \in E} a_{er}^k w_e, \forall r \in P^k, k = 1,\dots,K. \quad (9)$$

Therefore, $\theta(w)$ is computed by solving two sub-problems $AP_y$ and $AP_z$ stated as follows:

$$AP_y : \theta_y(w) = \text{Minimize} \sum_{e \in E} \sum_{l=1}^{L_e} \tilde{c}_e^l y_e^l \quad (10)$$

subject to: (2) and (6).

$$AP_z : \theta_z(w) = \text{Minimize} \sum_{k=1}^{K} \sum_{r=1}^{|P^k|} \tilde{h}_r^k z_r^k \quad (11)$$

subject to (3) and (5).

It is noteworthy that $AP_y$ has the intergrality property. Moreover, it can be solved by applying a simple inspection approach. In order to strengthen $AP_y$, we propose to amend it by adding the following valid inequalities:

$$\sum_{e=(o^k,j) \in E} \sum_{l=1}^{L_e} y_e^l \geq 1, \forall k = 1,\dots,K / o^k = s^k \text{ or } o^k = t^k. \quad (12)$$

Constraints (12) enforce that, from any source node or any sink node, at least one edge should be selected in a feasible solution.

Furthermore, we notice that solving $AP_z$ amounts to solve an all-pair shortest paths problem on the graph $G = (V, E)$ with modified arc costs $w_e$ for each edge $e$, $e \in E$. Since all these edge costs are nonnegative, the pricing problem can be efficiently solved using Floyd-Warshall all-pair shortest paths algorithm.

## 4 SOLVING THE LAGRANGIAN DUAL

To solve the Lagrangian dual problem, one of the most effective approaches is the well-know subgradient method. Its general statement is given below.

**Subgradient algorithm**

**Step 0**: *Initialization.* Start with an initial zero Lagrangian multipliers vector $w^0$ and $q=0$.

**Step 1**: *Subgradient computation.* Given $w^q$, solve the relaxed Lagrangian problem to obtain a solution $(y,z)$ and its corresponding cost $\theta^q$. Compute the current subgradient $g^q$ as follows:

$$g_e^q = \sum_{k=1}^{K} \sum_{r=1}^{|P^k|} a_{er}^k z_r^k - \sum_{l=1}^{L_e} u_e^l y_e^l, \forall e \in E. \quad (13)$$

If $g^q=0$ (practically if $||g^q||<10^{-6}$), then optimal solution is reached and stop. Otherwise, go to Step 2.

**Step 2**: *Search direction computation.* Select a search direction $d^q$ (as described in section 4.1). If $d^q=0$ (practically if $||d^q||<10^{-6}$), then put $d^q=g^q$ and go to Step 3.

**Step 3**: *Lagrangian multipliers computation.* Let $w^q=w^q+\lambda^q d^q$ for some step-lengh

$$\lambda^q = \beta^q \cdot \frac{UB - \theta^q}{\|d^q\|^2}, \quad (14)$$

where UB is an upper bound of the optimal solution and $\beta^q \in (0,2]$ is a step-length parameter. Set $q \leftarrow q+1$.

**Step 4**: *Termination test*. If a maximum number *MaxIt* of consecutive iterations is ran without any improvement of the Lagrangian dual function value then stop. Otherwise, go to Step 1.

In the following paragraphs, we describe different search direction and step-length strategies that we used to build the tested variants of the subgradient algorithm.

## 4.1 Search Direction Rules

In order to efficiently derive tight lower bounds for the DCMNDP, we explore six different rules to compute the search direction vector. These rules have been proposed by (Polyak, 1967, 1969; Camerini et al., 1975; Crowder, 1976; Sherali & Ulular, 1990; Junger et al., 1995) and the corresponding subgradient algorithms are denoted by **SG1-SG6**.

For all these variants, we initialize the search direction as $d^0 = g^0$. At each iteration $q$, in variants **SG1-SG5**, the search direction is updated as $d^q = g^q + \sigma^q d^{q-1}$ through the appropriate value of $\sigma$ detailed in Table 1. In variant **SG6**, the search direction is considered as a linear combination of the incumbent and previous subgradients and is computed as $d^q = \alpha\ g^q + (1-\alpha)g^{q-1}$ at each iteration $q>0$ (Junger et al., 1995). In our experiments, the parameter $\alpha$ is fixed to 0.7.

Table 1: Search direction computation rules.

| SG Algorithm | Value of $\sigma^q$ |
|---|---|
| **SG1** | *0*, (Polyak, 1967, 1969) |
| **SG2** | $\begin{cases} -1,5g^q d^{q-1}/\|d^{q-1}\|^2 & \text{if} \quad g^q d^{q-1} < 0 \\ 0 & \text{otherwise} \end{cases}$ (Camerini et al., 1975) |
| **SG3** | $\begin{cases} \|g^q\|/\|d^{q-1}\| & \text{if} \quad g^q d^{q-1} < 0 \\ 0 & \text{otherwise} \end{cases}$ (Camerini et al., 1975) |
| **SG4** | *0.8*, (Crowder, 1976) |
| **SG5** | $\|g^q\|/\|d^{q-1}\|$, (Sherali & Ulular, 1990) |

We notice here that **SG3** is commonly referred to as the modified-Camerini-Frata-Maffioli variant and **SG5** as the Average Direction Search (ADS) variant.

## 4.2 Step-length Parameter Rules

Since the convergence of the subgradient algorithm depends heavily on the step-length value $\lambda$, and more precisely the step-length parameter $\beta$, we investigated their impact on the computation time through testing the following three step-length parameter rules:

**R1** is the commonly used rule (Legendre & Minoux, 1977). It determines the coefficients $\beta^q$ dynamically taking into account the progression of the lower bound at each iteration $q$. First, we set $\beta^0=2$, then $\beta^q$ is halfed from its previous value whenever the Lagrangian dual function value fails to increase.

**R2** and **R3** were proposed in (Held & Karp, 1971), where $\beta^0=2$, then $\beta^q$ is halfed from its previous value every $p$ consecutive iterations. The value of $p$ depends on the size problem. In our case, we consider **R2** and **R3** with the value of $p$ fixed to *2n* and *2m*, respectively.

**R4**, **R5**, and **R6** were proposed by Polyak (Polyak, 1969), where $\beta^q$ is constant at each iteration $q$. We have tested several values of $\beta$. For the sake of conciseness, we will present the results of only 3 values of $\beta$: *0.01* (**R4**), *0.1* (**R5**), and *1.99* (**R6**).

## 5 COMPUTATIONAL RESULTS

The computational experiments aim to first, assess the empirical performance of the Lagrangian relaxation when applied to an arc-path formulation and second, compare different subgradient algorithms in terms of lower bound quality and computation time.

We implemented all the subgradient algorithms with Microsoft Visual Studio 2010 Ultimate C++ and we carried out all the computational experiments on a CORE i5 2.2 GHz Personal Computer with 8.0 GB RAM.

Experiments were conducted on a test-bed that consists of two sets of instances: 20 instances that are randomly generated according to the work of (Mrad & Haouari, 2008) and 12 real-world instances picked from the network design literature (Fumagalli et al., 1999; Miya & Saita, 1999; Walter 2002; Laland, 2004; *SNDLIB*). For these instances, the number of nodes and edges range from *10* to *64*, and *15* to *100*, respectively.

Tables 2-3 present the characteristics of the tested instances. The first column gives the name of the instance (**Inst.**), the second column provides the number of nodes (**n**), the third column presents the number of edges (**m**) and the forth column provides the number of commodities (**K**). For the randomly generated instances, the commodities are defined between each pair of nodes within the set $V$, i.e. $K=n(n-1)/2$ and $3$ facility types are available for each edge, i.e. $L_e = 3$, $\forall\ e \in E$. For the real world instances, the number of commodities ranges from 21 to 4032 and the number of facility types ranges from 2 to 40. The last column of Table 3 gives the number of facility types (**L**).

Table 2: Characteristics of the randomly generated instances.

| Inst. | n | m | K |
|---|---|---|---|
| D.1 | 10 | 15 | 45 |
| D.2 | 15 | 20 | 105 |
| D.3 | 15 | 25 | 105 |
| D.4 | 15 | 30 | 105 |
| D.5 | 20 | 35 | 190 |
| D.6 | 20 | 40 | 190 |
| D.7 | 20 | 45 | 190 |
| D.8 | 21 | 40 | 210 |
| D.9 | 22 | 45 | 231 |
| D.10 | 23 | 50 | 253 |
| D.11 | 24 | 55 | 276 |
| D.12 | 25 | 60 | 300 |
| D.13 | 25 | 50 | 300 |
| D.14 | 30 | 60 | 435 |
| D.15 | 34 | 70 | 595 |
| D.16 | 40 | 75 | 780 |
| D.17 | 40 | 75 | 780 |
| D.18 | 45 | 80 | 990 |
| D.19 | 50 | 90 | 1225 |
| D.20 | 50 | 100 | 1225 |

The set of real world problems comprises:
- Five instances denoted by (F1-F3), GRID35 and GRID12 provided by France Telecom (Lalande, 2004; Walter, 2002),
- One instance denoted by (NSF) and corresponds to the National Science Foundation Networks (Miya & Saita, 1999),
- One instance denoted by (EON) that is extracted from the European Optical Network (Fumagalli et al., 1999),
- One instance (EU) provided by (Lalande, 2004),
- Four instances (Pdh, Di-yuan, Nobel_us, Nobel_Germany) from *SNDLIB* (Orlowski., 2007).

Table 3: Characteristics of the real world instances.

| Inst. | n | m | K | L |
|---|---|---|---|---|
| F1 | 11 | 25 | 110 | 2 |
| Grid12 | 12 | 17 | 132 | 2 |
| NSF | 14 | 21 | 21 | 2 |
| EON | 18 | 37 | 37 | 2 |
| F2 | 20 | 34 | 34 | 2 |
| Grid35 | 35 | 58 | 1190 | 2 |
| F3 | 41 | 77 | 77 | 2 |
| EU | 64 | 81 | 4032 | 2 |
| Pdh | 11 | 34 | 24 | 3 |
| Di-Yuan | 11 | 42 | 22 | 6 |
| Nobel_us | 14 | 21 | 91 | 40 |
| Nobel_Geramany | 17 | 26 | 121 | 40 |

In our experiments, all tested subgradient algorithms are stopped after *100* consecutive iterations without any improvement of the Lagrangian dual function value (*MaxIt=100*).

In a first set of experiments, we compare the six variants of the subgradient algorithm **SG1-SG6**, and we use the popular Legendre & Minoux rule **R1** for computing the step-length parameter. The results are presented in Tables 4 and 5. We report the average GAP between the lower bound and the best one in percentage and the average CPU time in seconds, respectively, for the randomly generated instances (**Rand.**) and the real-world problems (**Real.**).

Table 4: Average GAP (%) using R1.

| Inst. | SG1 | SG2 | SG3 | SG4 | SG5 | SG6 |
|---|---|---|---|---|---|---|
| **Rand.** | 13.18 | 9.08 | **1.60** | 9.41 | 4.87 | 10.26 |
| **Real.** | 7.26 | 7.26 | 2.89 | 7.26 | **0.00** | 7.26 |
| **Aver.** | 10.22 | 8.17 | 2.24 | 8.33 | 2.43 | 8.76 |

Table 5: CPU time (sec) using R1.

| Inst. | SG1 | SG2 | SG3 | SG4 | SG5 | SG6 |
|---|---|---|---|---|---|---|
| **Rand.** | 0.86 | 0.91 | 1.12 | 0.88 | 1.33 | 0.85 |
| **Real.** | 1.48 | 1.34 | 1.18 | 1.46 | 1.33 | 1.23 |
| **Aver.** | 1.17 | 1.13 | 1.15 | 1.17 | 1.33 | 1.04 |

Tables 4 and 5 demonstrate the effectiveness of the Lagrangian relaxation applied to the arc-path formulation: lower bounds are derived in less than 1.33 seconds.

Furthermore, the results of table 4 show that when the step-length parameter **R1** is used, the best lower bound is obtained by the subgradient algorithm **SG3** for the randomly instances and by the subgradient algorithm **SG5** for the real world instances. Moreover, table 5 indicates that the CPU times are comparable for all the tested subgradient variants and they are around 1.20 seconds.

In a second set of experiments, we assess the impact of the different step-length parameter rules **R1-R6** on the lower bound quality and on the CPU time when the best two subgradient variants **SG3** and **SG5**, are used. Tables 6-9 present the average GAP between the lower bound and the best one in percentage and the average CPU time in seconds, respectively, for the randomly generated instances (**Rand.**) and the real-world problems (**Real.**).

Table 6: Average GAP (%) using SG3.

| Inst. | R1 | R2 | R3 | R4 | R5 | R6 |
|---|---|---|---|---|---|---|
| Rand. | 5.91 | 8.10 | 6.08 | **1.57** | 5.67 | 6.17 |
| Real. | 39.38 | 39.38 | 38.07 | **0.06** | 35.27 | 39.11 |
| Aver. | 22.65 | 23.74 | 22.08 | **0.82** | 20.47 | 22.64 |

Table 7: CPU time (sec) using SG3.

| Inst. | R1 | R2 | R3 | R4 | R5 | R6 |
|---|---|---|---|---|---|---|
| Rand. | 1.19 | 1.10 | 1.18 | 2.13 | 1.41 | 1.12 |
| Real. | 1.02 | 1.03 | 1.05 | 15.76 | 1.23 | 0.14 |
| Aver. | 1.10 | 1.07 | 1.11 | 8.94 | 1.32 | **0.63** |

Table 8: Average GAP (%) using SG5.

| Inst. | R1 | R2 | R3 | R4 | R5 | R6 |
|---|---|---|---|---|---|---|
| Rand. | 5.49 | 6.49 | 5.06 | **1.74** | 4.81 | 6.11 |
| Real. | 38.69 | 37.75 | 37.56 | **0.16** | 34.75 | 39.76 |
| Aver. | 22.09 | 22.12 | 21.31 | **0.95** | 19.78 | 22.93 |

Table 9: CPU time (sec) using SG5.

| Inst. | R1 | R2 | R3 | R4 | R5 | R6 |
|---|---|---|---|---|---|---|
| Rand. | 1.35 | 1.29 | 1.41 | 4.72 | 1.36 | 1.34 |
| Real. | 1.38 | 1.35 | 1.36 | 22.24 | 1.37 | 0.15 |
| Aver. | 1.36 | 1.32 | 1.38 | 13.48 | 1.36 | **0.74** |

Tables 6 and 8 show that when using subgradient algorithms **SG3** and **SG5**, the best lower bounds are obtained by using the step-length parameter rule **R4** with an average gap of 0.82% and 0.95%, respectively. Moreover, we observe that, the average gap increases with the value of $\beta$ for rules **R4-R6**.

However, according to tables 7 and 9, the rule **R4** requires the highest CPU times when compared to the other step-length parameter rules.

When comparing the two subgradient variants **SG3** and **SG5** along with the six step-length parameter rules **R1-R6**, tables 6-9 show that the best algorithm is **SG3** with rule **R4**. Indeed, it provides a good lower bound (within 0.82% of the best generated lower bound) within 8.94 seconds of CPU time.

An instructive finding is that the modified-Camerini-Frata-Maffioli variant (**SG3**) and the Average Direction Strategy (**SG5**) consistently outperformed all the other subgradient algorithms.

# 6 CONCLUSION

In this paper, we investigated the NP-hard Discrete Cost Multicommodity Network Design Problem (DCMNDP) which arises in several varieties of network design applications including telecommunications, cable television and logistic networks. Thus, solving such an NP-hard problem deserves considerable algorithmic challenges to researchers and requires the use of sophisticated strategies that aim to relax or/and decompose the problem and to strengthen its model with valid inequalities. In order to efficiently generate tight lower bounds, we applied the Lagrangian relaxation technique to an arc-path formulation for the DCMNDP, instead of the commonly used arc-node formulation. Then, to solve the obtained Lagrangian dual problem, we investigated the empirical performance of six deflected subgradient algorithms that are most suited for tackling large-scale dual Lagrangian problems since they require relatively few memory and computational effort. Indeed, we explored the deflected subgradient algorithms using various direction-search and step-length strategies, which provide a great deal of flexibility on determining lower bounds for the DCMNDP. We reported the results of extensive computational experiments carried out on randomly generated instances and on real-world problems. Thereby, the performance of the proposed deflected subgradient algorithms is compared. Future research will focus on embedding the proposed lower bound within a Branch-and-Bound scheme in order to derive exact solutions.

## ACKNOWLEDGEMENTS

The authors would express their gratitude to Professor Mohamed Haouari for his insightful suggestions.